\documentclass[graybox]{svmult}
\usepackage{geometry}
\usepackage{cite}
\usepackage{titlesec}
\titlespacing*{\section}{0pt}{2ex}{2ex}

\usepackage{amsmath,amssymb}
\DeclareMathOperator{\sign}{sign}
\usepackage{mathptmx}       % selects Times Roman as basic font
\usepackage{helvet}         % selects Helvetica as sans-serif font
\usepackage{courier}        % selects Courier as typewriter font
\usepackage{type1cm}        % activate if the above 3 fonts are
\usepackage{makeidx}         % allows index generation
\usepackage{graphicx}        % standard LaTeX graphics tool
\usepackage{multicol}        % used for the two-column index
\usepackage[bottom]{footmisc}% places footnotes at page bottom

\makeindex             % used for the subject index
\begin{document}
\renewcommand{\thesection}{9.\arabic{section}}
\renewcommand{\thedefinition}{9.\arabic{definition}}
\renewcommand{\theequation}{9.\arabic{equation}}
\renewcommand{\thefigure}{9.\arabic{figure}}
\renewcommand{\thetheorem}{9.\arabic{theorem}}
\renewcommand{\theproposition}{9.\arabic{proposition}}
\renewcommand{\thecorollary}{9.\arabic{corollary}}

\setcounter{page}{219}
\setcounter{chapter}{8}
{
\title{Chapter 9. \\ Theory of Differential Inclusions \\ and Its Application in Mechanics}
% Use \titlerunning{Short Title} for an abbreviated version of
% your contribution title if the original one is too long
%\let\clearpage\relax
}
\author{Maria Kiseleva, Nikolay Kuznetsov and Gennady Leonov}
% Use \authorrunning{Short Title} for an abbreviated version of
% your contribution title if the original one is too long
\institute{Maria Kiseleva \at St. Petersburg State University, Russia
\and Nikolay Kuznetsov \at St. Petersburg State University, Russia;
University of Jyv\"{a}skyl\"{a}, Finland
\\ email: nkuznetsov239@gmail.com
\and Gennady Leonov
\at Institute of Problems of Mechanical Engineering
Russian Academy of Sciences, \\
St.Petersburg, Russia
\and
\mbox{ \ } %\copyright \ Springer International Publishing AG 2018 
\\
J.B. Clement and W. Yu (eds.), {\it New Perspectives and Applications of Modern \\
Control Theory} \url{https://doi.org/10.1007/978-3-319-62464-8\_9}
}
% Use the package "url.sty" to avoid
% problems with special characters
% used in your e-mail or web address
%

\titlerunning{Theory of Differential Inclusions and Its Application in Mechanics}
\authorrunning{M. Kiseleva et al.}

\vspace{-1cm}
\maketitle
\vspace{-2cm}
\abstract*{The following chapter deals with systems of differential equations with
discontinuous right-hand sides. The key question is how to define the solutions of
such systems. The most adequate approach is to treat discontinuous systems as
systems with multivalued right-hand sides (differential inclusions). In this work
three well-known definitions of solution of discontinuous system are considered.
We will demonstrate the difference between these definitions and their
application to different mechanical problems. Mathematical models of drilling
systems with discontinuous friction torque characteristics are considered. Here,
opposite to classical Coulomb symmetric friction law, the friction torque
characteristic is asymmetrical. Problem of sudden load change is studied. Analytical methods of investigation of systems with such asymmetrical friction
based on the use of Lyapunov functions are demonstrated. The Watt governor and Chua system
are considered to show different aspects of computer modeling of
discontinuous systems.}

\abstract{The following chapter deals with systems of differential equations with
discontinuous right-hand sides. The key question is how to define the solutions of
such systems. The most adequate approach is to treat discontinuous systems as
systems with multivalued right-hand sides (differential inclusions). In this work
three well-known definitions of solution of discontinuous system are considered.
We will demonstrate the difference between these definitions and their
application to different mechanical problems. Mathematical models of drilling
systems with discontinuous friction torque characteristics are considered. Here,
opposite to classical Coulomb symmetric friction law, the friction torque
characteristic is asymmetrical. Problem of sudden load chance is studied. Analytical methods of investigation of systems with such asymmetrical friction
based on the use of Lyapunov functions are demonstrated. The Watt governor and Chua system
are considered to show different aspects of computer modeling of
discontinuous systems.}

\section{Introduction}
\begin{figure}
    \centering
        \includegraphics[width=0.5\textwidth]{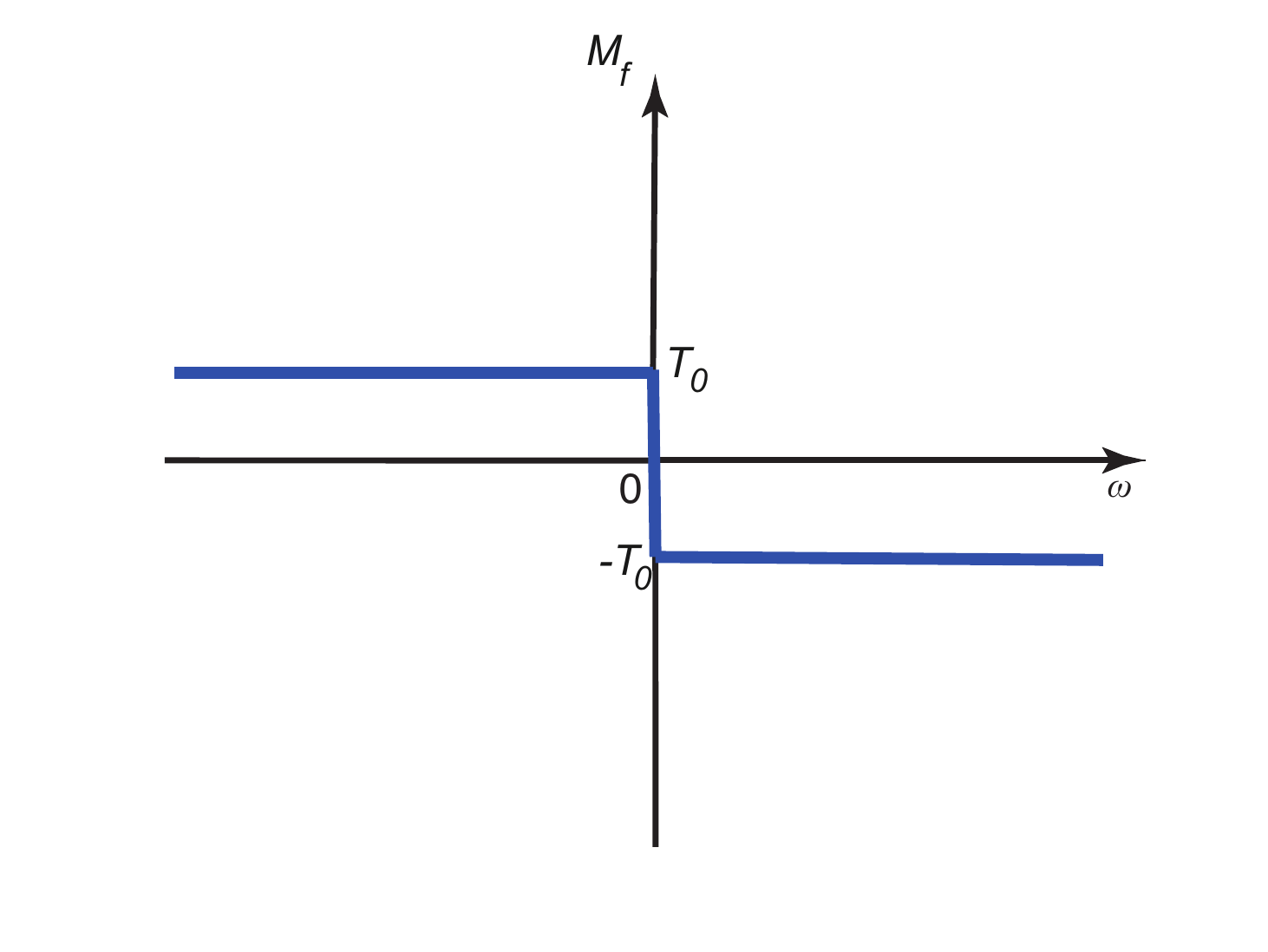}
    \caption{Coulomb friction torque $M_f$}
    \label{Classicalfriction}
\end{figure}

Two hundred and thirty years ago, after numerous experiments, Coulomb has formulated a law of dry fiction (Coulomb friction, see Fig. \ref{Classicalfriction}). Since then various problems stimulated the development of theory
of mechanical systems with dry friction.

First of all, it is important to mention the well-known Penleve paradoxes \cite{Penleve1954}, which provoked interesting discussions and showed contradiction of Coulomb's
\linebreak
law with Newton's laws of classical
mechanics. Nowadays, an independent research
\linebreak
branch named ``tribology'' has grown out of these
classical problems. Many
\linebreak
researches contributed to this branch of science, among them there are such
famous scientists as F.P. Bowden and D. Tabor \cite{Bowden1950}, E. Rabinovicz \cite{Rabinowicz1965}, P.J. Blau \cite{Blau2008}, K.C. Ludema \cite{Ludema1996}, I.G. Goryacheva \cite{Goryacheva1998,Goryacheva2001}, V.I. Kolesnikov \cite{Kolesnikov2003}.

From a mathematical point of view the problem of investigation of dynamics in models with dry
friction is closely connected with the theory of differential inclusions and dynamical systems with
discontinuous right-hand sides. Nowadays this theory
\linebreak
is being actively developed and applied to
investigation of different applications by such famous scientists as
S.V. Emelyanov, A.S. Poznyak, V.I. Utkin and others
(see e.g. \cite{Agrachev2004,Bennett1993,Brogliato1999,Pfeiffer1996,Utkin2013,Emelyanov1967,Plestan2010,Poznyak1999,Utkin2013,Edwards1998,Arkin1998,Kloeden2011}).
%theory of adaptive systems

\begin{figure}[t!]
    \centering
        \includegraphics[width=0.5\textwidth]{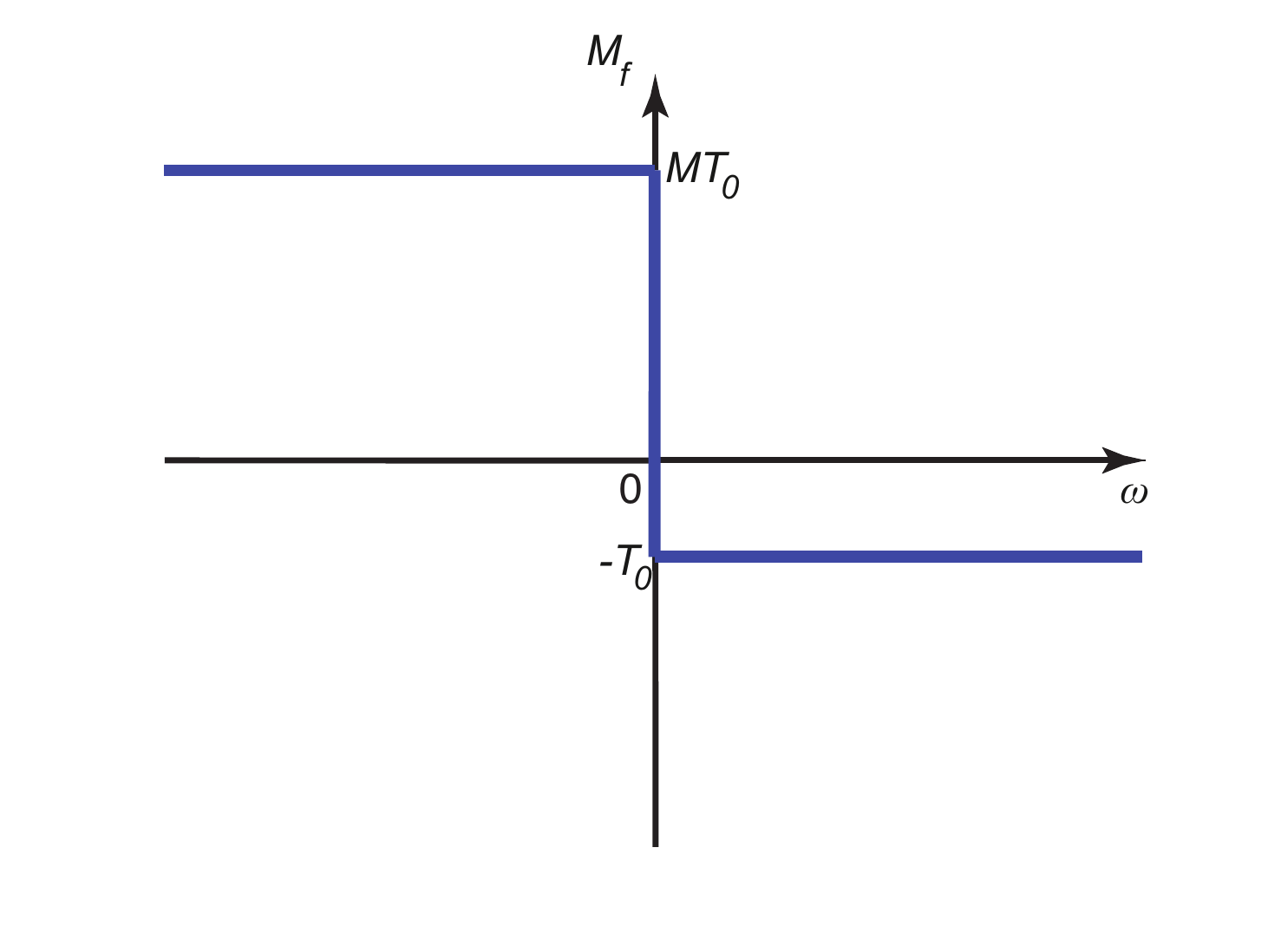}
    \caption{Non-symmetric friction torque $M_f$}
    \label{friction}
\end{figure}

The following work is motivated by the problem of investigation of a drilling system. This problem was studied by the research group from the Eindhoven Uni-\linebreak versity of Technology \cite{Mihailovic-2004,deBruin-2009}. In these papers, the interaction of the drill with the bedrock is described by symmetric discontinuous characteristics.
In the paper \cite{LeonovKiseleva-DAN-2012} a more precise model of friction is considered for simplified mathematical model of drilling system actuated by induction motor. Here the following assumption is made: the moment of resistance force with asymmetric characteristics (see Fig. \ref{friction}, $M$ is assumed to be large enough) is used instead of classical Coulomb friction with symmetric discontinuous characteristics. Such an asymmetric characteristic has a ``locking'' property -- it allows rotation of the drill in one direction only. The considered simplified model corresponds to an ordinary hand electric drill. In this case it \linebreak is
naturally to assume that the drilling takes place in one direction only.

The study of discontinuous systems with dry friction is a challenging task due to the need for a special theory for discontinuous systems to be developed. In particular, a proper definition of the solution on discontinuity surface is required. Now there are many definitions of solutions of discontinuous system, here three of them are consid\-ered following the works \cite{Gelig2006,GeligLY-1978,Filippov1985,KiselevaK-2015-Vestnik,LeonovKKK-2015-IFAC}. Analytical investigation of stability of simplified drilling systems will be performed. Additional examples of numeri\-cal modeling theory of differential inclusions and its application in discontinuous
{\parfillskip=0pt\par}

 \pagebreak\noindent
\noindent mechanical systems will be considered. It will be explained why it is necessary to use special methods of investigation for discontinuous systems.

\section{Differential Equation with Discontinuous Right-Hand Sides
and Differential Inclusions: Definitions of Solutions}

The starting point of studies in theory of differential inclusions is usually connected with the works of French mathematician A. Marchaud and Polish mathematician S.~K. Zaremba published in 1934-1936. They were studying equations of the form
\begin{equation}\label{diffincl1}
D x \subset f(t,x),
\end{equation}
where %$f: \mathbb{R} \times \mathbb{R}^n \to \mathbb{R}^n$ and
$t\in \mathcal{D}_{t} \subset \mathbb{R}$, $x \in \mathcal{D}_{x} \subset \mathbb{R}^n$  and $f(t,x)$
is a multivalued vector function that \linebreak maps each point $(t,x)$ of some region $\mathcal{D} = \mathcal{D}_{t} \times \mathcal{D}_{x}$  %from $\mathbb{R}^{n+1}$
to the set $f(t,x)$ of points from $\mathbb{R}^n$. For operator $D$ the notions of contingent and paratingent were introduced by Marchaud and Zaremba respectively.

\begin{definition}\label{defcont}
Contingent of vector function $x(t)$ at the point $t_0$
is a set Cont $x(t_0)$ of all limit points of sequences $\dfrac{x(t_i)-x(t_0)}{t_i-t_0}$, $t_i \rightarrow
t_0$, $i = 1,2,...$
\end{definition}

\begin{definition}\label{defparat}
Paratingent of vector function $x(t)$ at the point $t_0$
is a set Parat $x(t_0)$ of all limit points of sequences $\dfrac{x(t_i)-x(t_j)}{t_i-t_j}$, $t_i \rightarrow
t_0$, $t_j \rightarrow
t_0$, $i = 1,2,...$
\end{definition}

Wazhewski continued investigations of Marchaud and Zaremba and proved \cite{Wazewski1961}
that if $x(t)$ is a solution of differential inclusion \eqref{diffincl1}
in the sense of Marchaud then vector function $x(t)$ is \emph{absolutely continuous}.

{\parfillskip=0pt\par}

 \pagebreak\noindent

\begin{definition}
Let $I \subset \mathcal{D}_t \subset \mathbb{R}$ be an interval of time. Function $x(t):I \to \mathbb{R}^n$
is \emph{absolutely continuous} on $I$ if for every positive number $\varepsilon$ there is a positive number $\delta$ such that whenever a finite sequence of pairwise disjoint sub-intervals $(t_{1k}, t_{2k})$ of $I$ with $t_{1k}, t_{2k} \in I$ satisfies
\begin{equation*}
\sum_k (t_{2k}-t_{1k})<\delta
\end{equation*}
then
\begin{equation*}
\sum_k ||x(t_{2k})-x(t_{1k})||<\varepsilon.
\end{equation*}
\end{definition}

Important property of absolutely continuous function $x(t)$
is that $x(t)$
has derivative $\dot x(t)$
almost everywhere on $I$
(see, e.g. \cite{Royden1988}).
This property played a key role in the development of theory of differential inclusions and equations with discontinuous right-hand side since it allowed to avoid artificial constructions
in Definition~\ref{defcont} and \ref{defparat} and to consider usual derivative almost everywhere.

In 1960 paper \cite{Filippov1960} was published by A.~F. Filippov,
where he considered solutions of differential equations with discontinuous right-hand side
as absolutely continuous functions.
Filippov approach is one of the most popular among other notions of \linebreak solutions of
\emph{systems with discontinuous right-hand sides}.
%Let us start the discussion with consideration of this definition.
Following \cite{Filippov1960}, consider a system
\begin{equation}\label{sysdisc}
\dot x = f(t,x), \qquad t \in \mathbb{R}, x \in \mathbb{R}^n,
\end{equation}
where $f:\mathbb{R} \times \mathbb{R}^n \to \mathbb{R}^n$ is a piecewise continuous function such that measure of the set of discontinuity points is assumed to be zero.

\begin{definition}\label{filippovdef}
Vector function $x(t)$, defined on an interval $(t_1,t_2)$, is called a solution of (\ref{sysdisc})
if it is absolutely continuous and for almost all $t\in(t_1,t_2)$ vector $\dot x(t)$ is within minimal closed convex set, which contains all $f(t,x')$ when $x'$ is within almost all \linebreak
$\delta$-neighbourhood of the point $x(t)$ in $\mathbb{R}^n$ (for fixed $t$), i.e.
\begin{equation}
\dot x \in \prod_{\delta>0} \prod_{\mu N=0} {\rm conv} \, f(t,U(x(t),\delta) - N).
\end{equation}
\end{definition}

Consider the case when system (\ref{sysdisc}) is autonomous and vector function $f(x)$ is discontinuous on some smooth surface $S$ in $\mathbb{R}^n$ and continuous in the neighbourhood of this surface. Let there exist limits $f_+(x)$ and $f_-(x)$ of vector function $f(x)$ when a point $x$ approaches $S$ from one or another side.
Suppose that the vectors $f_+(x)$ and $f_-(x)$ are both pointing towards the discontinuity surface $S$. Then the so-called sliding mode appears. According to Definition \ref{filippovdef}, the vector field of sliding mode on the discontinuity surface can be defined as follows. The plane tangent to the surface $S$ at the point $x$ and the segment $l$, which connects the terminal points of vectors $f_+(x)$ and $f_-(x)$, are constructed. Then the vector with initial point at $x$ and terminal point at the point of intersection of the segment and tangent plane is constructed: $f_0=f_0(x)$.
According to Definition \ref{filippovdef}, vector $f_0(x)$ defines vector field at the point $x$.

\newpage

The obtained solution of (\ref{sysdisc}) satisfies Definition~\ref{filippovdef}, but nevertheless there are important applied problems for which Definition \ref{filippovdef} is unsuitable. As an example of such problem we consider a problem of synthesis of controls $u_1$ and $u_2$, which are limited, $|u_1| \leq 1$, $|u_2|\leq 1$, and which transform optimally fast each point $(x_1(0),x_2(0))$ of the system
\begin{equation}\label{syscontrol1}
\dot x_1 = x_2 u_1, \qquad \dot x_2 = u_2
\end{equation}
to the origin of coordinates. It is well-known \cite{Boltyansky-1969} that synthesis of such control is possible for the whole plane $(x_1,x_2)$. For example, for the first quadrant of the plane the optimal control is as follows
\begin{equation}
u_1 = \begin{cases}
1, &\text{ $x_1<0.5{x_2}^2$,}\\
-1, &\text{ $x_1\geq 0.5{x_2}^2$,}
\end{cases}
\qquad
u_2 = \begin{cases}
-1, &\text{ $x_1\leq 0.5{x_2}^2$,}\\
1, &\text{ $x_1> 0.5{x_2}^2$.}
\end{cases}
\end{equation}
In particular, the trajectory $x_1 = 0.5 {x_2}^2$ is optimal and for this trajectory system \linebreak(\ref{syscontrol1}) takes the form $\dot x_1 = -x_2$, $\dot x_2 = -1$.
Let us take the point $x = (x_1,x_2)$ on \linebreak this trajectory and approach to this trajectory from the side $x_1<0.5 {x_2}^2$. The limit value of the right-hand sides of system (\ref{syscontrol1}) is $f_+(x) = (x_2,-1)$. If we approach the trajectory from the side $x_1>0.5{x_2}^2$, then the limit is $f_-(x) = (-x_2,1)$. Since $f_+(x) = -f_-(x)$, in this particular case the segment $l$ passes through the point $x$,
\linebreak
i.e. $f_0(x) = 0$ and according to definition \ref{filippovdef} the solution on sliding mode is equilibrium state. At the same time $(-x_2,-1)$ is a velocity vector on optimal trajectory. Thus, optimal trajectory is not a solution in the sense of Definition \ref{filippovdef} by Filippov.

M.A. Aizerman and E.S. Pyatnitskiy \cite{AizermanP-1974}
offered other definition of solution
\linebreak
of equations with discontinuous right-hand sides which allows one to deal with
\linebreak
usual derivative. We consider their approach in the particular case when  $f(t,x)$ is discontinuous on the surface $\Sigma$.
Consider a sequence of continuous vector functions $f_\varepsilon (t,x)$, which coincide with $f(t,x)$ outside of $\varepsilon$-neighbourhood of surface $\Sigma$, and tend to $f(t,x)$ for $\varepsilon \to 0$ at each point, which does not belong to $\Sigma$. Let $x_\varepsilon (t)$ be a solution of the system
\begin{equation}
\dot x = f_\varepsilon(t,x).
\end{equation}
Then the solution of system (\ref{sysdisc}) in the sense of Aizerman and Pyatnitskiy is a limit of any uniformly converging subsequence of solutions $x_{\varepsilon_k}(t)$:
\begin{equation*}
x_{\varepsilon_k}(t) \rightrightarrows x(t).
\end{equation*}

In general, there may exist more than one such limit. Nevertheless this notion of solution, introduced in \cite{AizermanP-1974}, does not always suitable for applications.

\begin{figure}[h]
  \centering
  \includegraphics[scale=0.4]{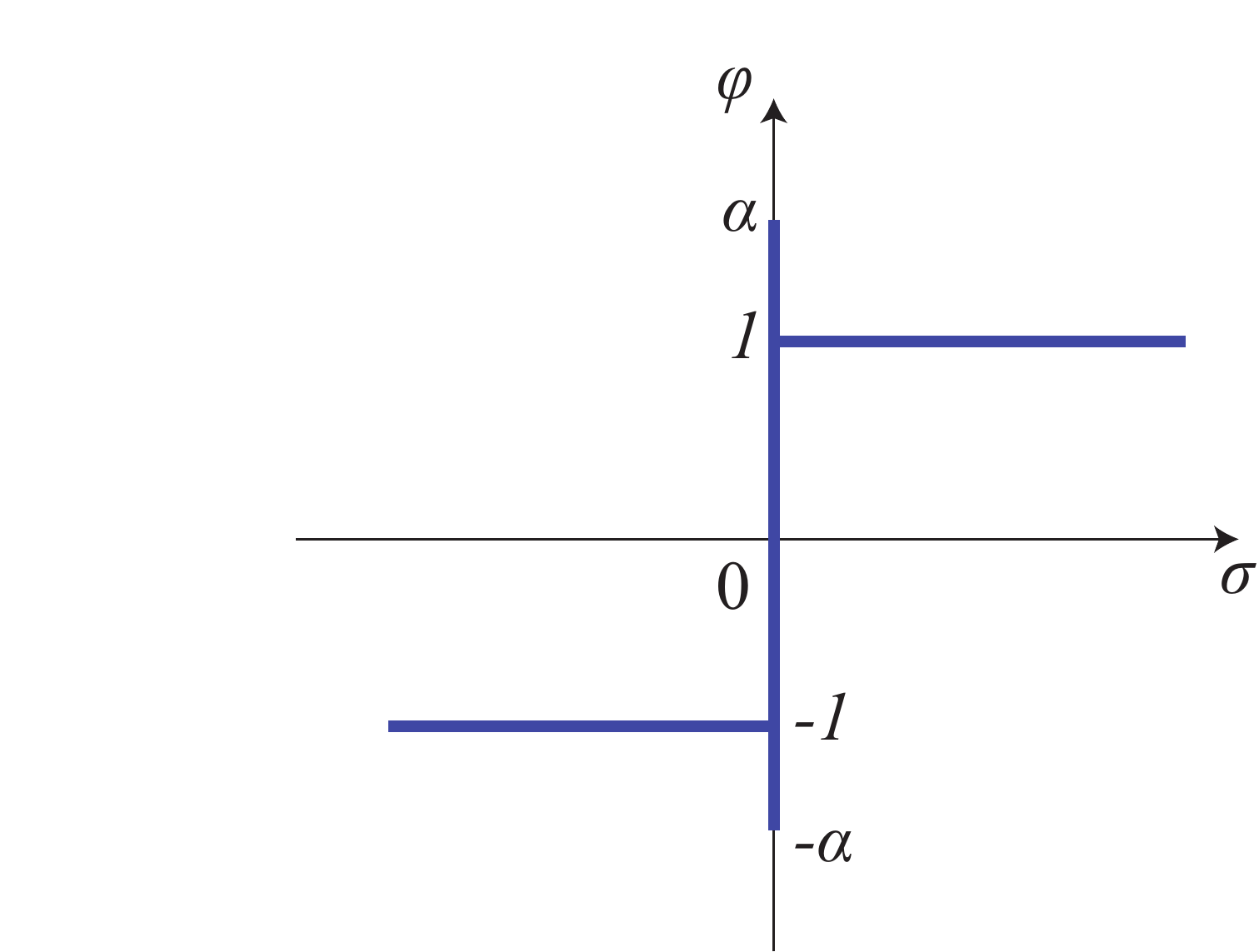}
  \caption{Friction model, where static friction exceeds dynamic friction}
  \label{example_dry}
\end{figure}

For example, consider a system
\begin{equation}\label{syscontrol}
\dot x = Ax + b \phi(\sigma), \quad \sigma = c^*x,
\end{equation}

{\parfillskip=0pt\par}

 \pagebreak\noindent
where $\phi(\sigma)$ is a dry friction characteristic,
shown in Fig.\ref{Classicalfriction} or in Fig.\ref{example_dry}, i.e.
\begin{equation}
\phi(\sigma) = \begin{cases}
\sign \sigma, &\text{ $\sigma \neq 0$,}\\
[-1,1], &\text{ $\sigma = 0$,}
\end{cases}
\quad \mbox{or} \quad
\phi(\sigma) = \begin{cases}
\sign \sigma, &\text{ $\sigma \neq 0$,}\\
[-\alpha,\alpha], &\text{ $\sigma = 0$.}
\end{cases}
\end{equation}

%and it is assumed that static friction values may be bigger then the values of dynamic friction, i.e. this is a friction with stalling characteristics
%. In this case the plot of the function $\phi(\sigma)$
%and it is of the form shown in Fig. \ref{example_dry}.

Since the definitions suggested by Filippov and by Aizerman and
Pyatnitskiy deal only with those values of a nonlinearity for which $\sigma \neq 0$, the solutions of system (\ref{syscontrol}) with dry friction characteristics,
shown in Fig.\ref{Classicalfriction} and in Fig.\ref{example_dry} coincide. This result does not match physics of this phenomena.

To take into account dynamics on the discontinuity surface,
the most adequate approach is to consider
system with discontinuous right-hand side \eqref{sysdisc}
as system with multivalued right-hand side, called \emph{differential inclusion} \cite{Filippov1971,GeligLY-1978}:
\begin{equation}\label{diffincl}
  \dot x \in f(t,x),
\end{equation}
where %$f: \mathbb{R} \times \mathbb{R}^n \to \mathbb{R}^n$ and
$t\in \mathcal{D}_{t} \subset \mathbb{R}$, $x \in \mathcal{D}_{x} \subset \mathbb{R}^n$
and $f(t,x)$ is a multivalued vector function that
\linebreak
maps each point $(t,x)$ of some region $\mathcal{D} = \mathcal{D}_{t} \times \mathcal{D}_{x}$  %from $\mathbb{R}^{n+1}$
to the set $f(t,x)$ of points from $\mathbb{R}^n$.
%where multivalued vector function $f(t,x)$ ($t\in \mathbb{R}^1$, $x \in \mathbb{R}^n$) takes each point $(t,x)$ of some region $\mathcal{D}$ from $\mathbb{R}^{n+1}$ to the set $f(t,x)$ of points from $\mathbb{R}^n$.
%The theory of differential inclusions was developed in the works .
%Let us bring some of the basic notions of this theory.

\begin{definition}\label{defDIsoution}
Vector function $x(t)$ is called a solution of differential inclusion
\linebreak
\eqref{diffincl},
if it is absolutely continuous
and for those $t$ for which derivative $\dot x(t)$ exists,
\linebreak
the following inclusion holds:
\begin{equation}\label{DIsoution}
  \dot x(t) \in f(t,x(t)).
\end{equation}
\end{definition}

%\begin{definition}
%Function $f(t,x)$ is called {\it continuous} at the point $(t_0,x_0)$ if it is semicontinuous and for any $\varepsilon>0$ there exists $\delta(\varepsilon,t_0,x_0)$ such that the set $f(t_0,x_0)$ is contained in $\varepsilon$-neighbourhood of set $f(t_1,x_1)$, provided that the point $(t_0,x_0)$ belongs to $\delta$-neighbourhood of point $(t_1,x_1)$.
%\end{definition}

%\section{Theory of differential inclusion}

To build a substantive theory it is assumed
that multivalued function $f(t,x)$ is semicontinuous.
\emph{Filippov approach} \cite{Filippov1985} requires additionally that $\forall (t,x) \in \mathcal{D}$
the set $f(t,x)$ is a minimal closed bounded set.
This conditions coincide with Definition~\ref{filippovdef}.

\begin{definition}
Function $f(t,x)$ is called  \emph{semicontinuous} (upper semicontinuous,
\linebreak
$\beta$-continuous) at the point $(t_0,x_0)$ if for any $\varepsilon>0$ there exists $\delta(\varepsilon,t,x)$ such that
{\parfillskip=0pt\par}

 \pagebreak\noindent
the set $f(t,x)$ is contained in the $\varepsilon$-neighbourhood of set $f(x_0,x_0)$, provided that
\linebreak
the point $(t,x)$ belongs to $\delta$-neighbourhood of the point $(t_0,x_0)$.
\end{definition}

As was shown above for some physical problems Filippov definition may gives wrong results, thus
a more general class of multivalued functions $f(t,x)$ was considered by
A.~Kh. Gelig, G.~A. Leonov and V.~A. Yakubovich \cite{GeligLY-1978}
(\emph{Gelig-Leonov-Yakubovich approach}): $\forall (t,x) \in \mathcal{D}$ the set
$f(t,x)$ is a bounded, closed, and convex set.

%We will refer to consideration of solution $x(t)$ of differential inclusion (\ref{diffincl}) with such condition on $f(t,x)$ as Gelig-Leonov-Yakubovich approach.
%Then solution $x(t)$ of differential inclusion (\ref{diffincl}) is called Gelig-Leonov-Yakubovich solution.

The following local theorem on the existence of solutions of differential inclusion holds true \cite{GeligLY-1978}.

\begin{theorem}
Suppose that multivalued function $f(t,x)$ is semicontinuous at
\linebreak
every point $(t_1,x_1)$ of a region
\begin{equation*}
\mathcal{D}_1 \subset \mathcal{D}: \qquad |t_1-t_0| \leq \alpha, \qquad |x_1-a| \leq \rho,
\end{equation*}
and set $f(t_1,x_1)$ is bounded, closed, and convex. In addition, suppose
\begin{equation*}
\sup |y| = c  \qquad \text{for} \qquad y \in f(t_1,x_1), \qquad (t_1,x_1) \in \mathcal{D}_1.
\end{equation*}
Then for $|t-t_0| \leq \tau = \min(\alpha,\rho/c)$
there exist at least one solution $x(t)$ with initial condition $x(t_0) = a$,
which satisfies (\ref{diffincl}) in the sense of Definition~\ref{defDIsoution}.
\end{theorem}

For differential inclusion (\ref{diffincl}) theorem on continuation of solution remaining
\linebreak
in
bounded region holds true. Also the theorem, which states that for every $\omega$-
\linebreak
limiting point of trajectory $\!x(t)\!$ there exists at least one trajectory that entirely consists
\linebreak
of $\omega$-limiting points and some other theorems of qualitative theory, are valid
\ \ \cite{GeligLY-1978,Filippov1985,Tolstonogov2012}.

For generalization of classic results of stability theory on solutions of differential inclusion (\ref{diffincl}) in Gelig-Leonov-Yakubovich approach the existence of procedure of determination of discontinuous right-hand side according to a chosen solution (i.e. existence of extended nonlinearity, which allows one to replace differential inclusion with differential equation) was proved by B.M.~Makarov specially for the monograph \cite{GeligLY-1978}.

Let us demonstrate now the methods of theory of differential equations with
\linebreak
discontinuous right-hand sides, described above, in concrete problems.

\section{Analytical methods of investigation of discontinuous systems: an example of mathematical model of drilling system with ``locking friction''}

Consider the simplified mathematical model of drilling system actuated by induction motor (here we follow the works \cite{LeonovKiseleva-DAN-2012,Kiseleva2013J}). Assume that the drill is absolutely rigid
{\parfillskip=0pt\par}

 \pagebreak\noindent
body stiffly connected to the rotor, which rotates by means of the magnetic field created by the stator of the induction motor. The value of interaction of the drill with the bedrock is defined as a value of resistance torque, which appears during the
\linebreak
drilling process. Such a system experiences rapidly changing loads during the drilling, thus it is necessary to investigate the behaviour of induction motor during load jumps, i.e. when resistance torque acting on the drill suddenly changes.

The following problem of stability is urgent since decrease of drilling systems failures plays important role in the oil and gas industry \cite{Vaisberg-2002,Shokir-2004}.
\linebreak
\mbox{\ \ \ } As the equations of electromechanical model of the drilling system~we~consider~the
\linebreak
equations of induction motor, proposed in \cite{LeonovK-2009,KondratevaLSR-2001}, supplemented with the resistance torque $M_f$ of drilling:

\begin{align}\label{asys}
\begin{split}
& L\dfrac{di_1}{dt} + Ri_1 = SB (sin \theta ) \dot \theta, \\
& L\dfrac{di_2}{dt} + Ri_2 = SB (cos \theta ) \dot \theta, \\
& I \ddot \theta \in -\beta SB (i_1 sin \theta + i_2 cos \theta) + M_f\bigg(\dfrac{R}{L}+\dot \theta\bigg).
\end{split}
 \end{align}

Here
$\theta$ is a rotation angle of the drill about the magnetic field created by the
\linebreak
stator, which rotates with a constant angular speed $\dfrac{R}{L}$, $i_1(t)$, $i_2(t)$ are currents in
\linebreak
rotor windings, $R$ is resistance of windings, $L$ is inductance of windings, $B$ is the induction of magnetic field, $S$ is an area of one wind, $I$ is an inertia torque of drill,
\linebreak
$\beta$ is a proportionality factor, $\omega = \dot{\theta}+\dfrac{R}{L}$ is an angular velocity of the drill rotation with respect to a fixed coordinate system. The resistance torque $M_f$ is assumed \linebreak
to be of the Coulomb type \cite{GeligLY-1978,Penleve1954}. Unlike the classic Coulomb friction law with symmetrical discontinuous characteristic the friction torque $M_f$ has non-symmetrical discontinuous characteristics shown in Fig. \ref{friction}.:

\begin{equation*}
M_f(\omega)=\begin{cases}
-T_0 &\text{for $\omega>0$}\\
[-T_0,MT_0] &\text{for $\omega=0$}\\
MT_0 &\text{for $\omega < 0$.}
\end{cases}
\end{equation*}

For $T_0\geq 0$ the number $M>$ is assumed to be large enough. That means that the drilling process only takes place when
$\omega>0$. Such characteristics does not allow for $\omega$ to switch from positive to negative values during the transient process in real drilling systems. In this case the system only gets stuck for $\omega = 0$ for a long enough period of time. These effects happen frequently during drilling operation and are
\linebreak
studied by the analysis of system (\ref{asys}).

Performing the nonsingular change of variables

 \pagebreak\noindent

  \begin{align*}
  & s = -\dot \theta,
  \\
  & x = \dfrac{L}{SB}(i_1 cos \theta - i_2 sin \theta),
  \\
  & y = \dfrac{L}{SB}(i_1 sin \theta + i_2 cos \theta)
 \end{align*}
we reduce system (\ref{asys}) to the following one:

\begin{equation}\label{sys}
  \begin{array}{l}
\dot s \in ay + \widetilde{M}_f(s),
\\
\dot y = -cy-s-xs,
\\
\dot x = -cx+ys,
\end{array} \end{equation}
where $\ a = \dfrac{\beta(SB)^2}{IL}, \ c = \dfrac{R}{L}$. Here variables $x$, $y$ define electric values in rotor wind\-ings and the variable $s$ defines the sliding of the rotor. $\widetilde{M}_f$ has the following form

\begin{equation*}
\widetilde{M}_f (s)=\begin{cases}
\gamma,&\text{ $s<c$;}\\
[-\gamma M,\gamma],&\text{ $s=c$;}\\
-\gamma M,&\text{ $s>c$;}\\
\end{cases}
\end{equation*}
where $\gamma = \dfrac{T_0}{I}$.

According to Makarov's theorem for any solution of (\ref{sys}) in the sense of the Gelig-Leonov-Yakubovich approach there exists extended nonlinearity $\widetilde{M}_{f0}$ such
\linebreak
that the following system is valid
\begin{equation}\label{sys1}
\begin{array}{l}
\dot s(t) = ay(t) + \widetilde{M}_{f0}(t),
\\
\dot y(t) = -cy(t)-s(t)-x(t)s(t),
\\
\dot x(t) = -cx(t)+y(t)s(t)
\end{array}
\end{equation}
for almost all $t$.

Let us conduct local analysis of equilibrium states of system (\ref{sys1}).

\begin{proposition} For $0\leq \gamma < \dfrac{a}{2}$ system (\ref{sys1}) has a unique asymptotically stable equilibrium state.
\end{proposition}

Indeed, for $\gamma = 0$ system (\ref{sys1}) has one asymptotically stable equilibrium
\linebreak
state $s = 0,$ $y = 0$, $x = 0$, which occurs when the rotation of drill with constant
\linebreak
angular speed is congruent to the rotation speed of the magnetic field (idle speed operation).

For $\gamma \in (0,\dfrac{a}{2})$ system (\ref{sys1}) has one equilibrium state
\begin{align*}
s_0 = \dfrac{c(a-\sqrt{a^2-4\gamma^2})}{2\gamma}, \qquad y_0 = -\dfrac{\gamma}{a}, \qquad x_0 = -\dfrac{\gamma s_0}{ac},
\end{align*}

 \pagebreak\noindent
where $s_0$ is the smallest root of the equation
\begin{align*}
\dfrac{acs}{c^2+s^2} = \gamma.
\end{align*}
In this case the direction of rotation of the drill and the magnetic field are the same, but the drill rotates with a lower angular speed $s_0<c$.

\bigskip

Assume that there is a sudden change in load at the moment $t = \tau$ from value $\gamma_0$ to value $\gamma_1$, where $0<\gamma_0<\gamma_1$. This occurs at the moment when the drill comes in contact with the bedrock. For $\gamma = \gamma_0$ the system experiences a unique state of stable equilibrium
\begin{equation*}
s_0 = \dfrac{c(a-\sqrt{a^2-4{\gamma_0}^2})}{2\gamma_0}, \quad y_0 = -\dfrac{\gamma_0}{a}, \quad x_0 = -\dfrac{\gamma_0 s_0}{ac}.
\end{equation*}

It is essential that after the transient process the solution $s(t)$, $x(t)$, $y(t)$ of
the system (\ref{sys}) with $\gamma = \gamma_1$ and the initial data $s(\tau) = \dfrac{c(a-\sqrt{a^2-4{\gamma_0}^2})}{2\gamma_0}$, $y(\tau) =$
$-\dfrac{\gamma_0}{a}$, $x(\tau) = -\dfrac{\gamma_0 s_0}{ac}$ tends to the equilibrium state
\begin{equation*}
s_1 = \dfrac{c(a-\sqrt{a^2-4{\gamma_1}^2})}{2\gamma_1}, \quad y_1 = -\dfrac{\gamma_1}{a}, \quad x_1 = -\dfrac{\gamma_1 s_1}{ac}
\end{equation*}
as $t\rightarrow+\infty$.

The following theorem holds.

\begin{theorem}
Let the following conditions be fulfilled
\begin{equation}\label{cond21}
\gamma_0<\dfrac{a}{2},
\end{equation}

\begin{equation}\label{cond2}
\gamma_1<min\left\{\dfrac{a}{2}, \quad 2c^2\right\},
\end{equation}

\begin{equation}\label{cond3}
3(M^2+2M)\gamma_1^2-8c^2\gamma_1 +3ac^2\geq0.
\end{equation}

Then the solution of system (\ref{sys1}) with $\gamma = \gamma_1$ and the initial data $s(\tau) = \dfrac{c(a-\sqrt{a^2-4{\gamma_0}^2})}{2\gamma_0}$, $y(\tau) = -\dfrac{\gamma_0}{a}$, $x(\tau) = -\dfrac{\gamma_0 s_0}{ac}$ tends to an equilibrium state of
\linebreak
this system as $t\rightarrow+\infty$.
\end{theorem}

Let us give the scheme of the proof of this theorem. We consider the region $\left\{s(t)<c\right\}$ of the phase space of system (\ref{sys1}).

Performing change of variables
 \pagebreak\noindent

$$
\eta(t) = ay(t) + \gamma_1, \qquad z(t) = -x(t)-\dfrac{\gamma_1}{ac}s(t),
$$
we reduce system (\ref{sys1}) to the following system
\begin{equation}\label{sys2}
  \begin{array}{l}
\dot s(t) = \eta(t), \\
\dot \eta(t) = -c \eta(t) + az(t)s(t) - \psi(s(t)), \qquad \text{ for almost all } t\\
\dot z(t) = -cz(t) - \dfrac{1}{a}s(t) \eta(t) - \dfrac{\gamma_1}{ac} \eta(t).
\end{array}
\end{equation}

Here $\psi(s) = -\dfrac{\gamma_1}{c}s^2+as -c\gamma_1$.

Introduce a function

\begin{equation*}
V(s,\eta,z) = \dfrac{a^2}{2}z^2+\dfrac{1}{2}\eta^2 + \int\limits_{s_1}^s \psi(s) ds.
\end{equation*}
For every solution of system (\ref{sys2}) from region $s(t)< c$ the following condition
\begin{equation}\label{diffV}
\dot V(s(t),\eta(t), z(t)) =
-a^2cz(t)^2-\dfrac{a\gamma_1}{c}\eta(t)z(t) - c\eta(t)^2 \leq 0
\quad \mbox{for almost all $t$ \ }
\end{equation}
is fulfilled.

Quadratic form in the right-hand side of system (\ref{diffV}) definitely negative under condition (\ref{cond2}).

We introduce a set

\begin{equation*}
\Omega = \left\{V(s,\eta,z)\leq \int\limits_{s_1}^c \psi(s) ds + \dfrac{(1+M)^2}{2}\gamma_1^2, s\in [s_2,c] \right\},
\end{equation*}
where the point $s_2<c$ is such that

\begin{equation*}
\int\limits_{s_2}^c \psi(s) ds + \dfrac{(1+M)^2}{2}\gamma_1^2 = 0.
\end{equation*}
The set $\Omega$ is limited and for $s = c$ it becomes

\begin{equation*}
\dfrac{a^2}{2}z^2+\dfrac{1}{2}\eta^2 \leq
\dfrac{(1+M)^2}{2}\gamma_1^2.
\end{equation*}
Returning to the initial coordinates ($x$, $y$, $s$), we obtain

\begin{equation*}
(x+\dfrac{\gamma_1}{a})^2+(y+\dfrac{\gamma_1}{a})^2 \leq \dfrac{(1+M)^2}{a^2}\gamma_1^2.
\end{equation*}

 \pagebreak\noindent
Note that this circle is below the upper boundary $y = \dfrac{M \gamma_1}{a}$ of the sliding region $\Delta = \left\{s = c, -\dfrac{\gamma_1}{a} \leq y \leq\dfrac{M \gamma_1}{a}\right\}$ of system (\ref{sys1}).

In the sliding region $\Delta$ system (\ref{sys1}) can be reduced to the system of ordinary differential equations

\begin{equation*}
  \begin{array}{l}
  \dot y(t) = -cy(t)-c-cx(t), \\
  \dot x(t) = -cx(t)+cy(t),

\end{array}
\end{equation*}
which is reduced by replacement of time $t = \dfrac{t_1}{c}$ to

\begin{equation}\label{syscs}
  \begin{array}{l}
  \dot{y}(t) = -y(t)-x(t)-1, \\
  \dot{x}(t) = -x(t)+y(t).
\end{array}
\end{equation}

We introduce a function

\begin{equation*}
W(x,y) = (x+\dfrac{\gamma_1}{a})^2+(y+\dfrac{\gamma_1}{a})^2.
\end{equation*}

Semicircles $\left\{W(x,y) =R^2, y>
-\dfrac{\gamma_1}{a}\right\}$, where $R\leq\dfrac{M+1}{a} \gamma_1$,
are noncontact
\linebreak
for system (\ref{syscs}). Indeed, for the solutions of the system (\ref{syscs}) under condition (\ref{cond2})
the following relation

\begin{equation*}
\begin{array}{l}
\dfrac{1}{2}\dot W(x(t),y(t)) = -y^2(t)-y(t)-\dfrac{\gamma_1}{a} + \dfrac{\gamma_1^2}{a^2} - (x(t)+\dfrac{\gamma_1}{a})^2
\\
= (\dfrac{2\gamma_1}{a}-1)y(t) +\dfrac{\gamma_1}{a}(\dfrac{2\gamma_1}{a}-1)-R^2<0
\end{array}
\end{equation*}
\bigskip
is valid.

The solution, which falls into the sliding region, necessarily comes out through the lower boundary $y = -\dfrac{\gamma_1}{a}$ into the region $s<c$ due to the fact that $\dot{s}<0$ for $s = c, y<-\dfrac{\gamma_1}{a}$. From condition (\ref{diffV}) it follows that this solution proves to be
\linebreak
inside the region $\left\{ V(s,\eta,z)\leq \int\limits_{s_1}^c \psi(s) ds\right\}$. Then it does not fall further into
the
\linebreak
sliding region, and tends to a unique equilibrium state $(s_1,y_1,x_1)$ of the system due to the limitation of $\Omega$. It is obvious that the trajectories, which fall into $\Omega$, but not existing in the sliding region, also tend to the equilibrium state.

 This allows one to prove that the system is dichotomic\footnote{System is called dichotomic if every solution bounded for $t>0$ tends to stationary set for $t \to +\infty$ \cite{LeonovK-2009,GeligLY-1978,Leonov-2001,Leonov-2006-ARC}} if condition (\ref{cond2}) is fulfilled (for details about classical results of Lyapunov see \cite{GeligLY-1978}).

 \pagebreak\noindent

The set $\Omega$ contains the point {$s = s_0$, $\eta = \gamma_1-\gamma_0$, $z = \dfrac{\gamma_0-\gamma_1}{ac}s_0$} if

\begin{equation}\label{ss}
\dfrac{(\gamma_1-\gamma_0)^2}{2c^2} {s_0}^2 +
\dfrac{(\gamma_1-\gamma_0)^2}{2} \leq \int\limits_{s_0}^c \psi(s) ds
+ \dfrac{(1+M)^2}{2}{\gamma_1}^2.
\end{equation}

For $\gamma_0<\gamma_1$ and condition (\ref{cond3}) we have

\begin{equation}\label{eq1}
\dfrac{(\gamma_1-\gamma_0)^2}{2} \leq \int\limits_{0}^c \psi(s) ds +
\dfrac{(1+M)^2}{2}{\gamma_1}^2.
\end{equation}

Let us show that

\begin{equation}\label{eq2}
\dfrac{(\gamma_1-\gamma_0)^2}{2c^2} {s_0}^2 \leq \int\limits_{s_0}^0 \psi(s) ds.
\end{equation}

Indeed, taking into account that $\gamma_0\leq 2\gamma_1$, we obtain:

\(
\dfrac{\gamma_1}{3c}{s_0}^2 - \dfrac{a}{2}s_0-\dfrac{(\gamma_1-\gamma_0)^2}{2c^2} s_0 + c \gamma_1 =
\dfrac{1}{12c^2 {\gamma_0}^2}(c^2(a-\sqrt{a^2-4{\gamma_0}^2})^2\gamma_1 - 3 a^2 \\
c^2 \gamma_0 + 3 a^2 c \sqrt{a^2-4{\gamma_0}^2} \gamma_0 - 3(\gamma_1-\gamma_0)^2 (a-\sqrt{a^2-4{\gamma_0}^2})\gamma_0 + 12 c^2 \gamma_1 {\gamma_0}^2) \geq \\
\dfrac{1}{12c^2 {\gamma_0}^2}(2 a^2 c^2 - 2 a c^2 \sqrt{a^2-4{\gamma_0}^2} \gamma_1 +3 a c^2\sqrt{a^2-4{\gamma_0}^2}- 3 a^2 c^2 \gamma_0 + \\
3 \sqrt{a^2-4{\gamma_0}^2} \gamma_0 {\gamma_1}^2 - 3 a {\gamma_1}^2 \gamma_0 + 8 c^2 {\gamma_0}^2 \gamma_1) \geq 0.
\)

Hence, from inequalities (\ref{eq1}) and (\ref{eq2}) we obtain condition (\ref{ss}).

Thus solution $s(t)$, $\eta(t)$, $z(t)$ with the initial data $s(\tau) = s_0$, $\eta(\tau) = \gamma_1-\gamma_0$, $z(\tau) = \dfrac{\gamma_0-\gamma_1}{ac}s_0$ tends to equilibrium state of the system.

Let $M$ be a reasonably large number such that condition (\ref{cond3}) of the theorem is fulfilled. In this case the following statement is valid.

\begin{corollary} Let $\gamma_0 = 0$ and
\begin{equation}\label{cond31}
\gamma_1<min\left\{\dfrac{a}{2}, \quad 2c^2\right\}.
\end{equation}
Then the solution of system (\ref{sys}) with $\gamma = \gamma_1$ and the initial data $s(\tau) = 0$, $y(\tau) = \\ 0, x(\tau) = 0$ tends to equilibrium state of this system as $t\rightarrow+\infty$.
\end{corollary}

\begin{figure}[ht]
    \centering
        \includegraphics[width=0.50\textwidth]{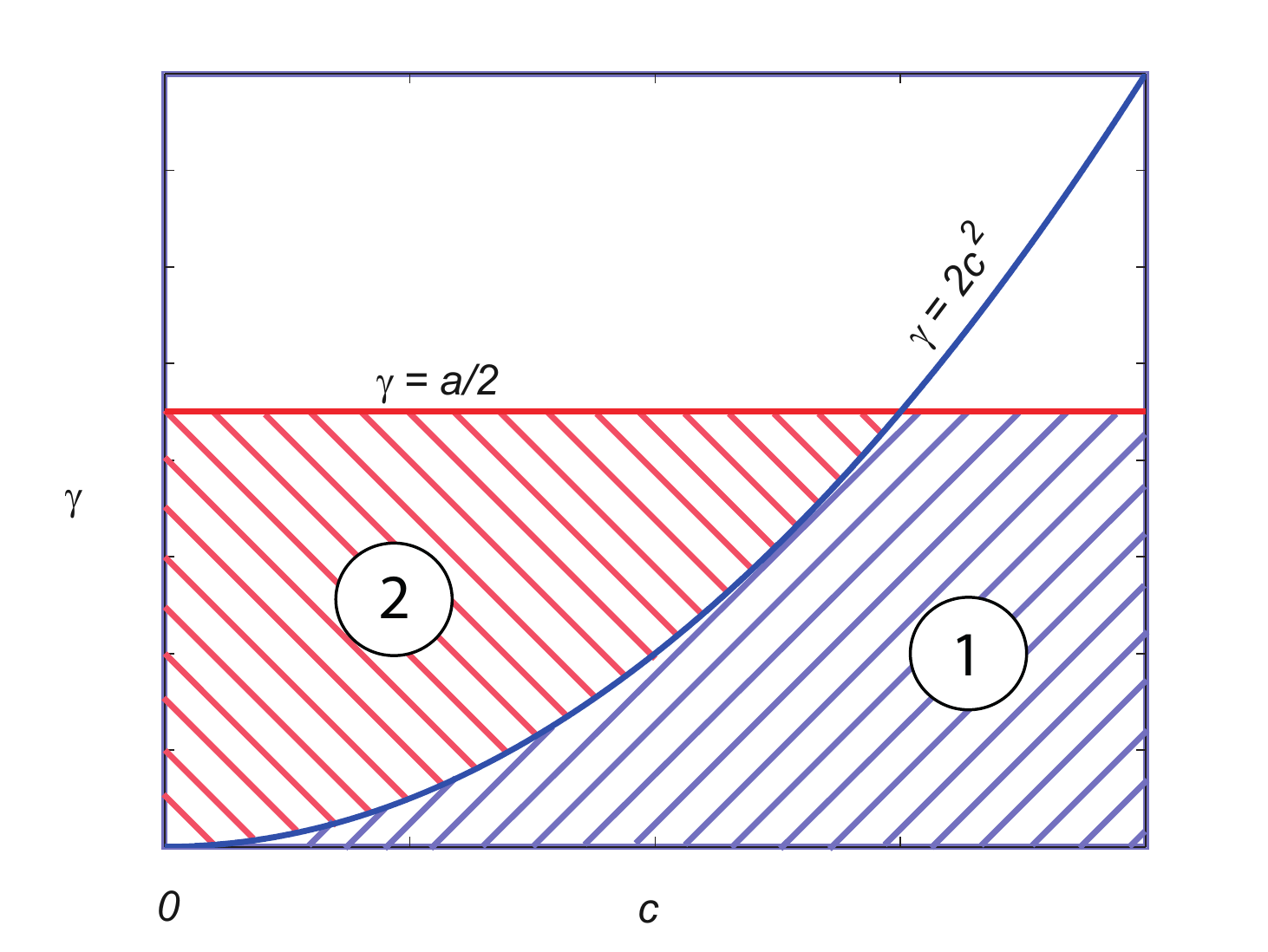}
    \caption{Safe load region: 1 -- due to the theorem, 2 -- obtained by numerical modeling of the system}
    \label{conditions}
\end{figure}

For the values $\gamma_1 \in \left\{ 2 c^2,\dfrac{a}{2}\right\}$ (i.e., condition (\ref{cond31}) is not fulfilled) the computer modeling of system (\ref{sys}) (region 2 in Fig. \ref{conditions}), which shows that the statement of consequence is retained, is carried out. Further we will discuss the aspects of modeling of systems with multivalued right-hand sides.

It can be checked that the extended nonlinearity can be written down in explicit form.

\pagebreak

\begin{corollary}
For system (\ref{sys}) the extended nonlinearity is of the following form:
\begin{equation*}
\widetilde{M}_{f0}=\begin{cases}
\gamma,&\text{if $s = c$, $y<-\dfrac{\gamma}{a}$ or $s<c$;}\\
-\gamma M,&\text{if $s = c$, $y>\dfrac{M \gamma}{a}$ or $s>c$;}\\
-ay,&\text{if $s = c$, $-\dfrac{\gamma}{a} \leq y \leq\dfrac{M \gamma}{a}$.}
\end{cases}
\end{equation*}
\end{corollary}
In the works \cite{KiselevaKL-2016-IFAC,KislevaKKL-2015-IFAC,KiselevaKKLS-2014,LeonovKKSZ-2014,KiselevaKLN-2014-DCNPS,LeonovKKN-2013-JAND} more complex models of drilling systems
\linebreak
were studied. Analytical investigation of such models is a challenging task, so it
\linebreak
is necessary to use
numerical methods. Let us further describe some aspects of
\linebreak
numerical modeling in two other applied systems -- Watt governor and Chua circuit.

\section{Numerical Methods of Investigation of Discontinuous \\ Systems}

\subsection{Difficulties of numerical modeling of discontinuous systems}
Numerical modeling is one of the tools of investigation of differential equations
\linebreak
with discontinuous equations
with right-hand sides.
Let us first show why it is
\linebreak
important to use special methods developed for discontinuous systems. Consider the
\linebreak
I.A. Vyshnegradsky problem.
The following system of differential equations
\linebreak
describes dynamics of Watt governor with dry friction

\begin{equation}\label{watt}
  \begin{aligned}
   & \dot y_1 = -Ay_1 + y_2 - \rm{sign}(y_1), \\
   & \dot y_2 = -B y_1 + y_3, \\
   & \dot y_3 = -y_1.
  \end{aligned}
\end{equation}
Let $\rm{sign}$ be understood here in ordinary sense:

\pagebreak
\begin{figure}
    \centering
        \includegraphics[width=0.70\textwidth]{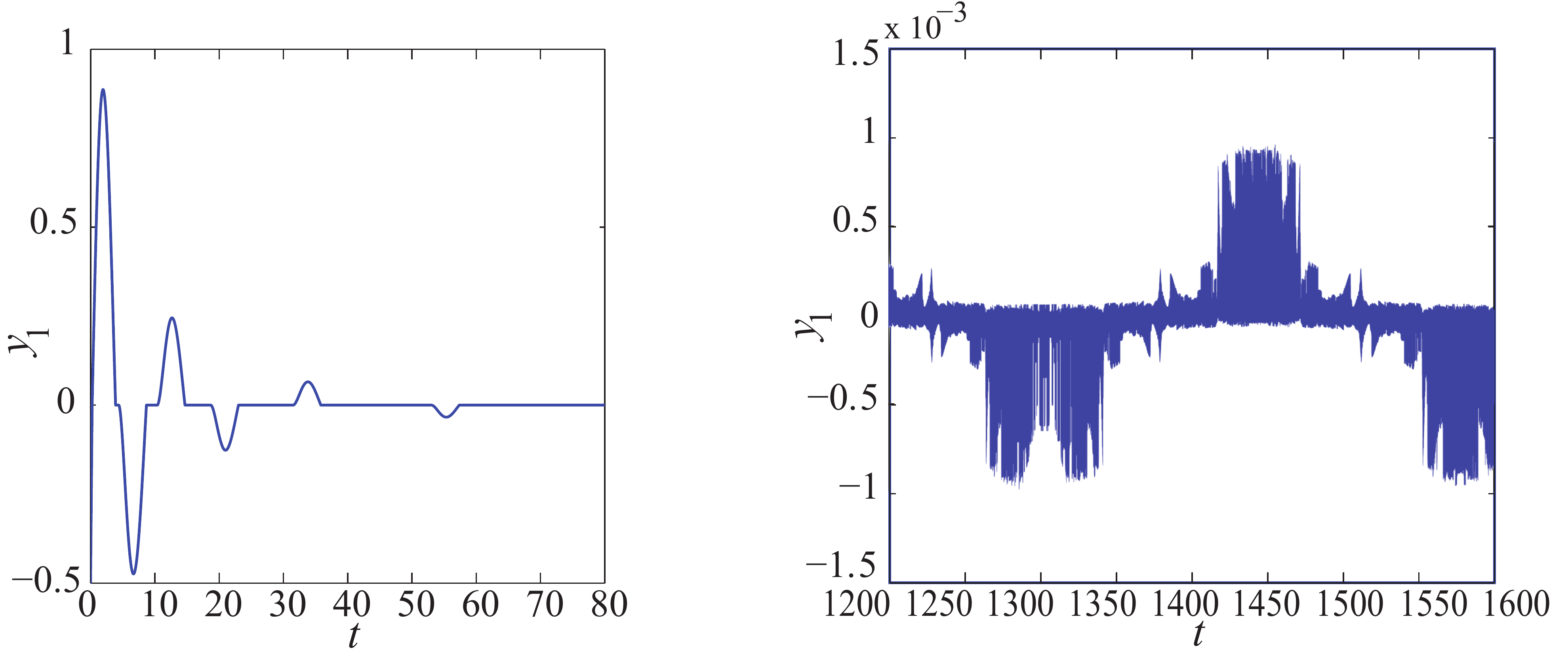}
    \caption{Numerical modeling of Watt governor using built-in Matlab functions}
    \label{WattError}
\end{figure}

\begin{equation*}
\rm{sign}(y_1)=\begin{cases}
1,&\text{if $y_1>0$;}\\
0,&\text{if $y_1=0$;}\\
-1,&\text{if $y_1<0$.}
\end{cases}
\end{equation*}
Let us consider the values of parameters $A = 1.5$, $B = 1.1$ and conduct numerical modeling of trajectory of system (\ref{watt}) with initial data $y_1(0)
= -0.5$, $y_2 = 1$, $y_3(0) = 1.2$, using standard Matlab build-in function {\it ode45} for solving ordinary differential equations. As can be seen in Fig. \ref{WattError} numerical modeling shows that there are oscillations in system (\ref{watt}).

This nonlinear system was studied by A.A. Andronov and A.G. Mayer \cite{AndronovM-1947}. In particular, they proved that sliding segment of this system is globally stable if the following inequalities
\begin{equation}
A>0, \qquad B>0, \qquad AB>1.
\end{equation}
are satisfied.

Thus, the result of modeling with standard build-in Matlab functions may lead to wrong results. Moreover, the notation (\ref{watt}) is wrong and right notation is as follows

\begin{equation}\label{watt1}
  \begin{aligned}
   & \dot y_1 \in -Ay_1 + y_2 - \rm{Sign}(y_1), \\
   & \dot y_2 = -B y_1 + y_3, \\
   & \dot y_3 = -y_1,
  \end{aligned}
\end{equation}
and the model of dry friction is described in the following way
\begin{equation*}
\rm{Sign}(y_1)=\begin{cases}
1,&\text{if $y_1>0$;}\\
[-1,1],&\text{if $y_1=0$;}\\
-1,&\text{if $y_1<0$.}
\end{cases}
\end{equation*}

\pagebreak
\begin{figure}
    \centering
        \includegraphics[width=0.70\textwidth]{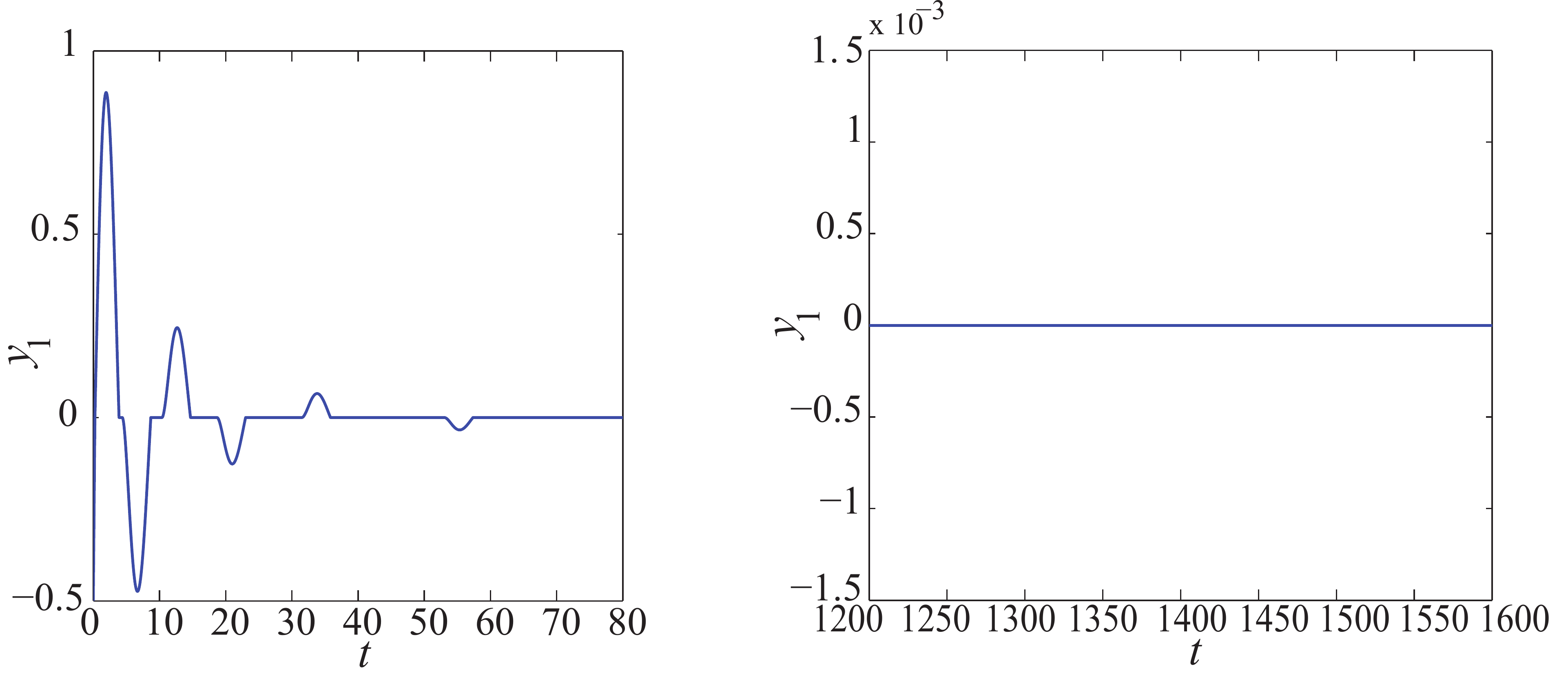}
    \caption{Numerical modeling of Watt governor using Filippov definition}
    \label{WattFilippov}
\end{figure}
\begin{figure}[h]
  \centering
  \includegraphics[scale=0.18]{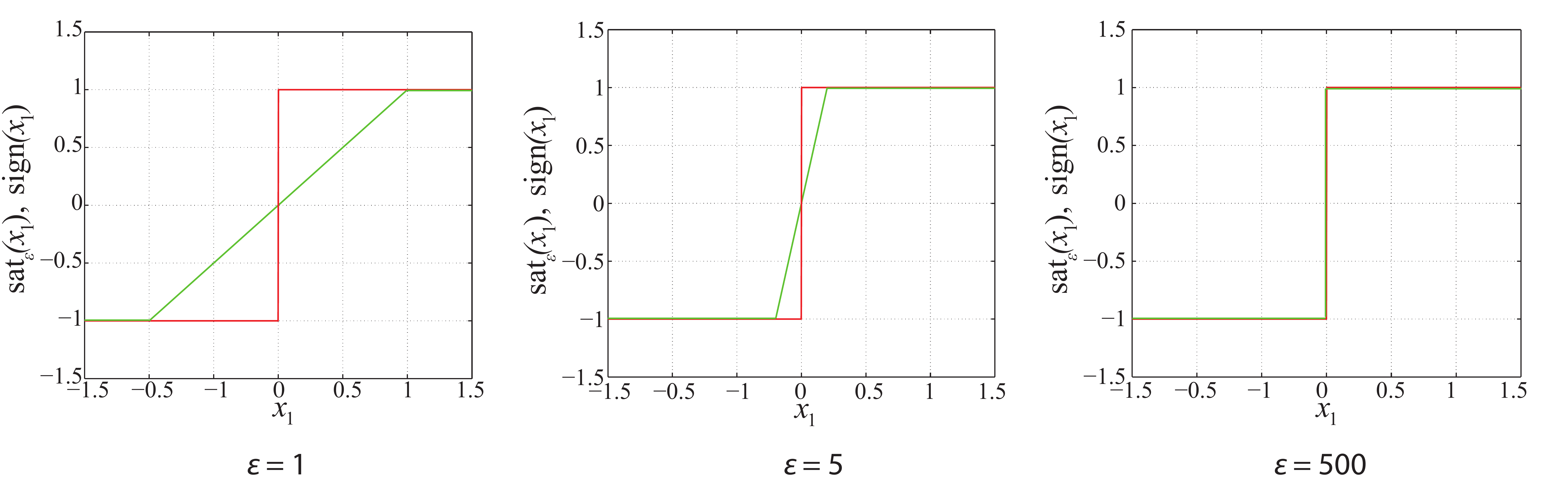}
  \caption{Plots of ${\rm sat}_\varepsilon(x_1)$ and ${\rm sign} (x_1)$}
  \label{example_SignSatChua}
\end{figure}

Let us conduct the numerical modeling using Filippov definition \cite{piiroinen2008event}.  The results of the modeling of system (\ref{watt1}) correspond t\texttt{}o theoretical results and are shown in Fig. \ref{WattFilippov}.

\subsection{Numerical Modeling of Chua System}

We showed an example of numerical modeling of a discontinuous system based
\linebreak
on Filippov approach.
Let us compare this method with modeling based on the Aizerman-Pyatnitskiy approach.

Consider the following example of discontinuous system --
modified Chua system with discontinues characteristic \cite{KuznetsovKLV-2013,KuznetsovKLV-2011-ICINCO,LeonovKKK-2015-IFAC}
\begin{equation}\label{chua}
\begin{array}{l}
\dot x_1 \in -\alpha(m_1 + 1)x_1 + \alpha x_2 - \alpha(m_0-m_1) {\rm Sign}(x_1), \\
\dot x_2 = x_1 - x_2 + x_3, \\
\dot x_3 = -\beta x_2  - \gamma x_3,
\end{array}
\end{equation}
where $\alpha$, $\beta$, $\gamma$, $m_0$, $m_1$ are parameters of the system.

\pagebreak
\begin{figure}[h]
  \centering
  \includegraphics[width=\textwidth]{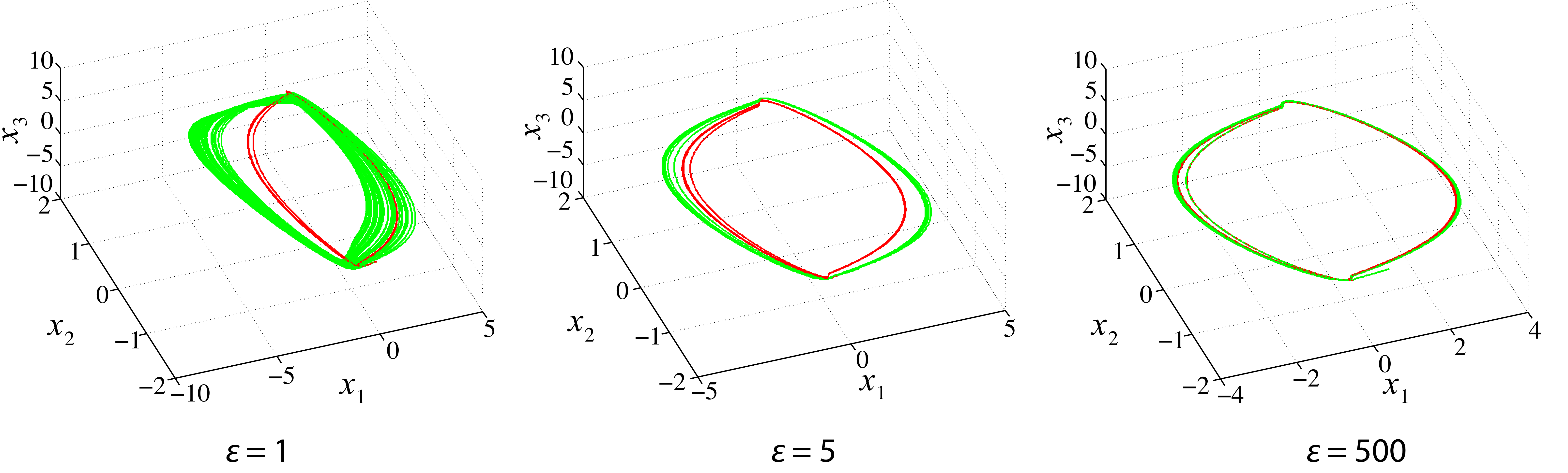}
 \caption{Hidden attractor in Chua system: comparison of solutions}
  \label{example_chua}
\end{figure}

For parameters $\alpha = 8.4562$, $\beta = 12.0732$, $\gamma = 0.0052$, $m_0 = -0.1768$, $m_1 = -1.1468$ system (\ref{chua}) has a so-called
{\it hidden attractor}\footnote{\!From a computational point of view, it is natural
to suggest the following classification of attractors,
based on the simplicity of finding the basin of attraction in the phase space:
\cite{KuznetsovLV-2010-IFAC,LeonovKV-2011-PLA,LeonovKV-2012-PhysD,LeonovK-2013-IJBC,LeonovKM-2015-EPJST,Kuznetsov-2016}:
{\it
 An attractor is called a \emph{self-excited attractor}
 if its basin of attraction
 intersects with any open
 \linebreak
 neighbourhood of an equilibrium,
 otherwise it is called a \emph{hidden attractor}.
}\!For a \emph{self-excited attractor}
\linebreak
its basin of attraction
is connected with an unstable equilibrium
and, therefore, self-excited attractors
\linebreak
can be localized numerically by the
\emph{standard computational procedure}
in which after a transient
\linebreak
process a trajectory,
started in a neighbourhood of an unstable equilibrium,
is attracted to the state of
\linebreak
oscillation and traces it.
For a hidden attractor its basin of attraction is not connected with equilibria.
\linebreak
The hidden attractors, for example, are the attractors in the systems
with no equilibria or with only one stable equilibrium
(a special case of multistability -- multistable systems and coexistence of attractors)
\cite{KuznetsovL-2014-IFACWC,DudkowskiJKKLP-2016,KiselevaKL-2016-IFAC}.
The classical examples of the hidden oscillations
are nested limit cycles in the study of 16th Hilbert problem
(see, e.g., \cite{KuznetsovKL-2013-DEDS}),
counterexamples to the Aizerman and Kalman conjectures on the absolute stability of nonlinear control systems \cite{LeonovBK-2010-DAN,LeonovK-2011-DAN,BraginVKL-2011},
and oscillations in electromechanical models without equilibria and with Sommerfeld effect
\cite{KiselevaKL-2016-IFAC}.
Recent examples of
\linebreak
hidden \mbox{\ } attractors \mbox{\ } can \mbox{\ } be \mbox{\ } found,\mbox{\ } e.g., \mbox{\ } in \mbox{\ }
\cite{AndrievskyKLS-2013-IFAC,AndrievskyKLP-2013-IFAC,KiselevaKKLS-2014,PhamJVWG-2014-HA,PhamRFF-2014-HA,WeiWL-2014-HA,LiSprott-2014-HA,ZhusubaliyevM-2015-HA,WeiML-2014-HA,LeonovKKSZ-2014,KuznetsovLYY-2014-IFAC,KuznetsovKMS-2015-HA,BurkinK-2014-HA,LiZY-2014-HA,ZhaoLD-2014-HA,LaoSJS-2014-HA,ChaudhuriP-2014-HA,LeonovKM-2015-CNSNS,MenacerLC-2016-HA,LiZCWH-2016-HA,PhamVJVK-2016-HA,Danca-2016-HA,JafariPGMK-2016-HA,Zelinka-2016-HA,LeonovKM-2015-EPJST}.
} \cite{KuznetsovLV-2010-IFAC,LeonovKV-2011-PLA,LeonovKV-2012-PhysD,LeonovK-2013-IJBC,LeonovKM-2015-EPJST,Kuznetsov-2016}.

So as to model system (\ref{chua})
with the help of both Filippov and Gelig-Leonov-Yakubovich definitions,
the special event-driven numerical method,
described in \cite{piiroinen2008event}, was used.
For modeling the system, using Aizerman-Pyatnitskiy approach, one needs to replace ${\rm sign}(x_1)$ by
${\rm sat}_\varepsilon(x_1) = \dfrac{1}{2}\bigg(|\dfrac{x_1}{\varepsilon}+1| - |\dfrac{x_1}{\varepsilon}-1|\bigg)$,
where $\varepsilon >0$. Decrease of parameter $\varepsilon$ allows one
to obtain Aizerman-Pyatnitskiy solution
(${\rm sat}_\varepsilon(x_1) \rightrightarrows {\rm sign}(x_1)$
for $\varepsilon \neq 0$, see Fig. \ref{example_SignSatChua}).
In Fig.\ref{example_chua} hidden attractor modeled
using Filippov definition method is drawn in red colour and hidden attractor modeled using Aizerman and Pyatnitskiy definition method is drawn in green.
As one can see,
the more $\varepsilon$ is decreased, the more solutions (attractors)
coincide with each other. This fact meets theorem proved in \cite{LeonovKKK-2015-IFAC,KiselevaK-2015-Vestnik}.

\vspace{-2mm}
\section*{\fontsize{11pt}{7.2}\selectfont Conclusion}
\vspace{-2mm}
We have discussed Filippov, Aizerman-Pyatnitskiy, and Gelig-Leonov-Yakubovich approaches
to the study of differential equation with discontinuous right-hand sides and differential inclusions. While for a wide range of dynamical models with these three approaches give the same result (see, e.g. \cite{LeonovKKK-2015-IFAC,KiselevaK-2015-Vestnik}),
there are models,\,where the

{\parfillskip=0pt\par}

 \pagebreak\noindent
difference between these definitions is essential.As examples,\,we have considered
the
\linebreak
Chua circuit with discontinues characteristic,\,and mechanical systems with classical
\linebreak
Coulomb symmetric friction law
and the asymmetrical friction torque characteristic.

\bigskip

\renewcommand{\abstractname}{Acknowledgements}
{\fontsize{8.5pt}{7.2}\selectfont
\begin{abstract}
 {This work was supported by the Russian Science Foundation (project 14-21-00041)}
\end{abstract}
}
%\bibliographystyle{spmpsci}
%\bibliography{C:/Dropbox/bib/bib-discont,C:/Dropbox/bib/bib_nk,C:/Dropbox/bib/bib_leonov,C:/Dropbox/bib/bib_full,C:/Dropbox/bib/bib-hidden}
%\bibliography{bib_KiselevaKL-2017-Springer,bib_nk,bib_full}

\end{document}